\documentclass[11pt,twoside]{amsart}
\usepackage{amsmath, amsthm, amscd, amsfonts, amssymb, graphicx, color}
\usepackage[bookmarksnumbered, colorlinks, plainpages]{hyperref}

\newtheorem{thm}{Theorem}[section]

\newtheorem{lem}[thm]{Lemma}
\newtheorem{prop}[thm]{Proposition}

\newtheorem{rem}[thm]{\bf{Remark}}

\newtheorem{example}[thm]{\bf{Example}}
\numberwithin{equation}{section}
\def\pn{\par\noindent}

\begin{document}


\vspace{1.3 cm}

\title{On Silverman's conjecture for a family of elliptic curves}
\author{Farzali Izadi and Kamran Nabardi}

\thanks{{\scriptsize
\hskip -0.4 true cm MSC(2010): Primary: 11G05; Secondary:  14H52.
\newline Keywords: Silverman's Conjecture, Elliptic Curve, Quadratic Twist, Rank, Parity Conjecture.\\
}}
 \maketitle
\begin{abstract}
 Let  $E$ be an elliptic curve over $\Bbb{Q}$ with the given Weierstrass equation  $ y^2=x^3+ax+b$. If $D$ is a squarefree integer, then let $E^{(D)}$ denote the $D$-quadratic twist of $E$ that is given by
   $E^{(D)}: y^2=x^3+aD^2x+bD^3$. Let $E^{(D)}(\Bbb{Q})$ be the group of $\Bbb{Q}$-rational points of $E^{(D)}$.
 It is conjectured by J. Silverman that there are infinitely many primes $p$ for which
  $E^{(p)}(\Bbb{Q})$ has positive rank, and there are infinitely
 many primes $q$ for which $E^{(q)}(\Bbb{Q})$ has rank $0$. In this paper, assuming the parity conjecture,
  we show that for infinitely many primes $p$, the elliptic curve $E_n^{(p)}: y^2=x^3-np^2x$ has
   odd rank and for infinitely many primes $p$, $E_n^{(p)}(\Bbb{Q})$ has even rank,
   where $n$ is a positive integer that can be written as biquadrates  sums in two different ways,
     i.e., $n=u^4+v^4=r^4+s^4$, where $u, v, r, s$ are positive integers such that $\gcd(u,v)=\gcd(r,s)=1$. More precisely, we prove that: if $n$  can be written in two different ways as biquartic sums and $p$ is prime, then under the assumption of the parity conjecture $E_n^{(p)}(\Bbb{Q})$ has odd rank  (and so a positive rank) as long as $n$ is odd and $p\equiv5, 7\pmod{8}$ or $n$ is even and $p\equiv1\pmod{4}$.
      In the end, we also compute the ranks of some specific values of $n$ and $p$ explicitly.
\end{abstract}
 \vskip 0.2 true cm
\pagestyle{myheadings}
\markboth{\rightline {\scriptsize  Izadi and Nabardi}}
        {\leftline{\scriptsize On Silverman's conjecture for a family of elliptic curves}}
\bigskip
\bigskip

\section{\bf Introduction}
\vskip 0.4 true cm An elliptic curve $E$ over  the rational field
$\Bbb{Q}$ is the projective curve associated to an affine equation
of the form
\begin{equation}
y^2=x^3+ax+b, \quad a,b\in \Bbb{Q},
\end{equation}
 where the cubic polynomial $x^3+ax+b$ has distinct roots. This is called the Weierastrass normal form.
 By the Mordell-Weil theorem it is known that the set of rational points  $E(\Bbb{Q})$ is a
finitely generated abelian group, as
such it takes the following decomposition:
\begin{equation}
E(\Bbb{Q})\simeq \Bbb{Z}^r \oplus  E(\Bbb{Q})_{{\rm{tors}}},
\end{equation}
where $r$ is a nonnegative integer called the
 $ rank$ of $E$ and $E(\Bbb{Q})_{{\rm{tors}}}$ is the finite abelian group consisting
 of all the elements of finite order in $E(\Bbb{Q})$ (  See \cite[Theorem 6.7, page 239]{Silv} for more information).\\
Let us  recall the definition of the  quadratic twist of an elliptic curve.
If $D$ is a squarefree integer, then let $E^{(D)}$ denote the $D$-quadratic twist of $E$ that is given by
   $E^{(D)}: y^2=x^3+aD^2x+bD^3$. The group of $\Bbb{Q}$-rational points of $E^{(D)}$ is shown by $E^{(D)}(\Bbb{Q})$.
\\

In \cite{IKN}, the authors  considered the family of elliptic curves
defined by
 \begin{equation}\label{E1}
E_n: y^2=x^3-nx,
\end{equation}
 where $n$ is a positive integer that can be written as biquadrates sums in two different ways,
  i.e.,  $n=u^4+v^4=r^4+s^4$ where $\gcd(u,v)=\gcd(r,s)=1$. The 4-tuple  $(u, v, r, s)$ with positive integer coordinates satisfying the above conditions
   is called a primitive solution.
 In  \cite[Theorem 1.1]{IKN} the authors proved  that  for such  $n$ the elliptic curve  $y^2=x^3-nx$, has rank $\ge 3.$ If moreover $n$ is odd and the parity conjecture is true, then it has even rank $\geq4$.\\

 The Diophantine equation
 \begin{equation}\label{E11}
 n=u^4+v^4=r^4+s^4,
 \end{equation}
 was first studied by Euler \cite{Eu1} in 1772 and since then  has been considered by numerous mathematicians. Among quartic Diophantine equations, \eqref{E11} has a distinct  feature for its simple structure, the almost perfect symmetry between the variables and close relationship with the theory of elliptic functions. The simplest parametric solution of \eqref{E11} that was constructed by Euler
  \cite[Page 201, (13.7.11)]{Ha} is:
 \begin{equation}\label{E12}
 \left\{\begin{array}{l}
 u=a^7+a^5b^2-2a^3b^4+3a^2b^5+ab^6,\\
 v= a^6b-3a^5b^2-2a^4b^3+a^2b^5+b^7,\\
 r=a^7+a^5b^2-2a^3b^4+3a^2b^5+ab^6,\\
 s=a^6b+3a^5b^2-2a^4b^3+a^2b^5+b^7.\end{array}\right.
 \end{equation}
 The first known examples of solutions, and among these the solution in least positive integers, i.e., $(u,v, r, s)=(134, 133, 158, 59)$, were cmputed already by Euler \cite{Eu1, Eu2, Eu3}. But this parametric solution does not exhaust   all the possible
solutions of the equation \eqref{E11} ( See \cite[page 21]{Hardy}). Some others were found by later researchers ( See \cite[pp. 644-647]{Dic} ), but it was not until the advent of computers that systematic searches could be conducted. The most extensive lists published to date are due to Lander and Parkin \cite{Lander1} and Lander, Parkin and Selfridge \cite{Lander2}.  Zajta \cite{Zajta}, discusses the  method for finding such  solutions  and presents a list of 218 numerical solutions in  the range $\max (u, v, r, s)<10^6$.
  Choudhry \cite{Chou},  presents a method of deriving new solutions of equation \eqref{E11} starting from a given  solution. According  to his method,  by taking $(\pm133,\pm134,\pm158,\pm59)$ one  can obtain:
\begin{equation*}
\begin{array}{cccc}
u&v&r&s\\
1054067&545991&522059&1057167\\
10381&10203&12231&2903\\
1453319&829418&1486969&461882\\
1137493&654854&60779&1167518\\
114613&111637&134413&34813\\
6565526&3687711&6710751&1967986\\
12178821457&7038985479&783453421&12505169907
\end{array}
\end{equation*}

\section{\bf {\bf \em{\bf Silverman's  conjecture and  some known results related to it}}}
\vskip 0.4 true cm

In this section, we first state  Silverman's conjecture and then we briefly discuss some known results related to it.\\

\textbf{\em Silveman's  Conjecture} (\cite[page250]{Ono1}).
\emph{If $E$ is an elliptic curve over the rational field $\Bbb{Q},$
then there are infinitely many primes $p$ for which $E^{(p)}(\Bbb{Q}),$ has
positive rank,
and there are infinitely many primes $q$ for which $E^{(q)}(\Bbb{Q})$ has  rank $0$} .\\

 By using $2$-descents ( See \cite[ Theorem 3.1, page 229]{Se}), one can prove part of this conjecture for the congruent number elliptic curve
 \begin{equation*}
 E': y^2=x^3-x.
 \end{equation*}
 For instance it is known for prime $p$  that if $p\equiv3\pmod{8}$, then $E'^{(p)}$ has rank $0$ and if $p\equiv5\pmod{8}$,
  then $E'^{(2p)}(\Bbb{Q})$ has rank $0$. Moreover, it is known that $E'^{(pqr)}(\Bbb{Q})$ has rank $0$ if $p$, $q$ and $r$ are primes satisfying
 \begin{equation*}
  p\equiv1\pmod{8},\  q\equiv3\pmod{8},\ r\equiv 3\pmod{8},\  \text{and}\  \left(\frac{p}{q}\right)=-\left(\frac{p}{r}\right),
  \end{equation*}
  where $\left(\frac{p}{q}\right)$ is the Legendre symbol.\\

  On the other hand, for the second part of the conjecture, i.e., the part with positive rank, Monsky \cite[Corollary 5.15, page 66]{Mon}
  proved that if $p_1$, $p_3$, $p_5$  and $p_7$ denote primes $\equiv1,\ 3,\ 5,\ 7\pmod{8}$,
  then the $D$-quadratic twist of $E^{'}$ has positive rank, where $D$ runs through  the following numbers:
  \begin{align*}
  &p_5,\ p_7,\ 2p_7\ \text{and}\  2p_3,
  \end{align*}
  \begin{align*}
  &p_3p_7,\ p_3p_5,\ 2p_3p_5\ \text{and}\ 2p_5p_7,\label{E999}
  \end{align*}
  \begin{align*}
  &p_1p_5\ \text{provided}\left(\frac{p_1}{p_5}\right)=-1,\ \  p_1p_7\ \text{and}\ 2p_1p_7\ \text{provided}\left(\frac{p_1}{p_7}\right)=-1,
  \end{align*}
  and
  \begin{align*}
  &2p_1p_3\ \text{provided}\ \left(\frac{p_1}{p_3}\right)=-1.
  \end{align*}
  Ono \cite[ Corollary 3.1, page 349]{Ono1},
  showed  that for some special curves $E$ there is a set $S$ of primes
   $p$ with density $\frac{1}{3}$ for which if $D=\prod p_j$ is a squarefree integer where $p_j\in S$,
    then $E^{(D)}$ has rank $0$. In particular $E^{(p)}$ has rank $0$ for every $p\in S$.

\section{\bf {\bf \em{\bf Main results}}}
\vskip 0.4 true cm Our main results are the followings:
\begin{thm}\label{t1}
Let $n$ be an odd number of the form \eqref{E11}. Then under the
assumption of the parity conjecture we have the following:
\begin{itemize}
\item[(i)]For any primes $p\equiv5,\ 7\pmod{8}$, the rank of $E_n^{(p)}(\Bbb{Q})$ is  odd and therefore
positive. Consequently, the  first part of Silverman's conjecture is true.
\item[(ii)]For any primes $p\equiv 1,\ 3\pmod{8}$, the rank of $E_n^{(p)}(\Bbb{Q})$ is even.
\end{itemize}
\end{thm}
\begin{thm}\label{t2}
Let $n$ be an even number of the form \eqref{E11}. Then under the
assumption of the parity conjecture we have the following:
\begin{itemize}
\item[(i)] For any primes  $p\equiv1\pmod{4}$, the rank of $E_n^{(p)}(\Bbb{Q})$ is odd and therefore
positive. Consequently, the  first part of Silverman's conjecture is true.
\item[(ii)] For any primes $p\equiv3\pmod{4}$, the rank of  $E_n^{(p)}(\Bbb{Q})$ is even.
\end{itemize}
\end{thm}
Before giving the proof of these theorems, we state a couple of
necessary facts from the literature.


 First of all,  the {\it Parity Conjecture} which states that an elliptic curve $E$ over $\Bbb{Q}$ with
  the  rank $r$
 satisfies
 \begin{equation*}
 \omega(E)=(-1)^r,
 \end{equation*} where $\omega(E)$ is the sign of the functional equation of the Hasse-Weil $L$-function
  $L(E,s)$.
  It is known  \cite{Birch} that for $n \not\equiv 0\pmod{4},$ a fourth power free integer, the sign of the functional equation,
  denoted $\omega(E_n),$ for the elliptic curve
   $$ E_n: y^2=x^3-nx,$$  is given by
\begin{equation}\label{E3.11}
\omega(E_n)={\rm{sgn}}(-n)\cdot\epsilon(n)\cdot\prod_{l^2\|n}\left(\frac{-1}{l}\right),
\end{equation}
where  $l\geq 3$  denotes a prime and

\begin{equation}\label{E3.1}
\epsilon(n)=\left\{
  \begin{array}{ll}
   -1, & n \equiv 1,3,11,13  \pmod{16}, \\
   &\\
    1, & n \equiv 2,5,6,7,9,10,14,15 \pmod{16}.\\
  \end{array}
\right.
\end{equation}
\vspace{.2cm}
(see Ono and Ono \cite[Equations $(1)$, $(2)$, and $(3)$]{Ono2}).\\
\hspace*{0.5cm} Second, the following proposition and  Theorem
\ref{3}  are also useful in the proof of our main results. Let us recall the lemma which  is essential to prove  Proposition \ref{2}.
\begin{lem}\label{l1}
Let $n$ be a nonzero integer, and let $p$ be an odd prime not dividing $n$. Then
\begin{equation*}
p\mid x^2+ny^2,\ \gcd(x,y)=1\Longleftrightarrow \left(\frac{-n}{p}\right)=1.
\end{equation*}
\end{lem}
\begin{proof}
 See (\cite [ Lemma 1.7, page 13]{Co}).
\end{proof}
\begin{prop}\label{2}
Let $n=u^4+v^4=r^4+s^4$ be such that \mbox{$\gcd(u,v)=1$}. If $l\ |\
n$ for an odd prime number $l$, then $l=8k+1$ for some
$k\in\Bbb{Z}$.
\end{prop}
\begin{proof}
Without loss of generality we can assume that $n$ is not divisible
by $4$. We use Lemma \ref{l1}.
Let $l$ be an odd prime factor of $n$. One can write
 \begin{equation*}
 n=u^4+v^4=(u^2-v^2)^2+2(uv)^2,
 \end{equation*}
 so,
   \begin{equation*}
   l\ |\ (u^2-v^2)^2+2(uv)^2.
   \end{equation*}
    According to Lemma \ref{l1}, $x=u^2-v^2$, $y=uv$ and $m=2$. Therefore, $\left(\frac{-2}{l}\right)=1$ which implies that $l=8k+1$ or $l=8k+3$. On the other hand, $n=(u^2+v^2)^2-2(uv)^2$, so
    \begin{equation*}
    l\ |\ (u^2+v^2)^2-2(uv)^2.
    \end{equation*}
   We get $\left(\frac{2}{l}\right)=1$ which implies that $l=8k+1$ or $l=8k+7$. Putting these two results together we get $l=8k+1$.
\end{proof}
\begin{rem}\label{rem1}
 By Proposition \ref{2}, every prime number $l$ dividing $n$ is in the form  $8k+1$ which in turn is in the form  $4k+1$ as well.
  Therefore by the quadratic reciprocity law ( See \cite[page153]{Silv1} ),
   for every prime factor $l$ of $n$ such that $l^2\|n$ we have $\left(\frac{-1}{l}\right)=1$.
   This shows that these prime factors can be ignored in evaluation of $\omega(E_n)$.
\end{rem}
\begin{thm}\label{3}
(Dirichlet's Theorem) If $n$ is a  positive integer  and $a$ and $b$
have no common divisor except $1$, then there are infinitely many
primes of the form $an+b$.
\end{thm}
\begin{proof}
See (\cite[Theorem 15, page 13]{Ha}).
\end{proof}
By the above facts in our disposal, we are now ready to investigate
Silverman's  conjecture  for the elliptic curve \eqref{E1}. To this
end, we first take the $p$-quadratic twists
\begin{equation}\label{E111}
E_n^{(p)}: y^2=x^3-p^2nx,
\end{equation}
where $n$ is  in the form  \eqref{E11} and the primes $p$ satisfying  $\gcd(p,n)=1,$ and then compute the sign of $\omega(E_n^{(p)})$ to characterize the parity of the rank for each curve.\\

It is a well-known fact that every odd prime  $p$ can be represented
as $p\equiv1\pmod{4}$, or $p\equiv3\pmod{4}$. Furthermore one can
easily check that square of every prime  $p$ can be written as
$p^2\equiv1\pmod{16}$, or $p^2\equiv9\pmod{16}$. It is also clear
that for each $p$ and each $n$ we have ${\rm {sgn}}(-p^2n)=-1$.
Having said that, we are now ready to prove our results.\\

{\bf {\bf \em{\bf Proof of Theorem 3.1 }}}\begin{itemize}
\item[(i)]
In this case, $u$ and $v$ have opposite parities. Without loss of
generality, let $u\equiv1\pmod{2}$ and $v\equiv0\pmod{2}$, then we
have
\begin{equation*}
n=u^4+v^4\equiv1\pmod{16}.
\end{equation*}
Next, based on the different choices for the primes  $p$, we get the
different results for $\omega(E_n^{(p)})$.
We have the following possibilities:\\
\\
(a):  $p\equiv1\pmod{4}$ and $p^2\equiv9\pmod{16}$.\\

Since $n\equiv1\pmod{16}$, we have $p^2n\equiv9\pmod{16}$ and then
from \eqref{E3.1} we have  $\epsilon(p^2n)=1$ implying that
$\omega(E_n^{(p)})=-1$ from \eqref{E3.11}. Therefore, under the assumption of the parity conjecture the rank is odd
and indeed it is  positive.
 For primes of the form  $p\equiv1\pmod{4}$, we have   $p=8k+1$, or $p=8k+5$ for some $k\in\Bbb{Z}$.
 These facts along with $p^2\equiv9\pmod{16}$ implies that
 $p\equiv5\pmod{8}$. By Dirichlet's theorem, there are infinitely many such primes and then, the first part of  Silverman's
 conjecture is true.\\

  (b): $p\equiv3\pmod{4}$ and $p^2\equiv1\pmod{16}$.\\

Similarly for this case from \eqref{E3.1}, we get
${\displaystyle\epsilon}(p^2n)=-1$. Moreover,
$\left(\frac{-1}{p}\right)=-1$.
 So  by \eqref{E3.11}, $\omega(E_n^{(p)})=-1$ meaning that the rank is odd (under the assumption of the parity conjecture) and indeed it is positive again.
 If $p\equiv3\pmod{4}$, then $p=8k+3$ or $p=8k+7$, for some $k\in\Bbb{Z}$.
 From these two and $p^2\equiv1\pmod{16}$, we get $p=8k+7$ and by Dirichlet's theorem, there are infinitely
 many such primes. Therefore, in this case the  first part of  Silverman's conjecture is true as well.
 \\
 \item[(ii)]
Like the previous case, we have
\begin{equation*}
n\equiv1\pmod{16},
\end{equation*}
and so, we can consider  two different cases as follows:\\

(a): {$p\equiv1\pmod{4}$ and $p^2\equiv1\pmod{16}$.}\\

Since $n\equiv1\pmod{16}$, so $p^2n\equiv1\pmod{16}$. Thus, \eqref{E3.1}
implies  that $\displaystyle{\epsilon}(p^2n)=-1$. It is clear that
$\left(\frac{-1}{p}\right)=1$.
In this case  from \eqref{E3.11} we have, $\omega(E_n^{(p)})=1$.
 Therefore, the parity conjecture states  the rank of the elliptic curve \eqref{E111} is  even.\\

(b): $p\equiv3\pmod{4}$ and $p^2\equiv9\pmod{16}$.\\

In this case one can easily check  that  $\epsilon(p^2n)$ in
\eqref{E3.1} equals to $1$. Moreover,
$\left(\frac{-1}{p}\right)=-1$,
 and so from \eqref{E3.11},  $\omega(E_n^{(p)})=1$. Therefore the rank is even as well.\\
 \end{itemize}

{\bf {\bf \em{\bf Proof of Theorem 3.2}}}
\begin{itemize}
\item[(i)]
 In this case, $u$ and $v$ are both odd and then
$n=u^4+v^4\equiv2\pmod{16}$. Now, as the previous theorem, we
consider different cases for the primes $p$ as  follows:\\

(a):  $p\equiv1\pmod{4}$ and $p^2\equiv9\pmod{16}$.\\

It is clear that $p^2n\equiv2\pmod{16}$, and so  by \eqref{E3.1},
$\epsilon(p^2n)=1$. Therefore, \eqref{E3.11} shows  that
$\omega(E_n^{(p)})=-1$ and so, under the assumption of  the parity conjecture the rank is positive.
 We have seen that these primes are in the form of $8k+5$ and there are infinitely many such primes.
It turns out the first part of   Silverman's conjecture is true.\\

(b):  $p\equiv1\pmod{4}$ and $p^2\equiv1\pmod{16}$.\\

Obviously, $p^2n\equiv2\pmod{16}$  and so from \eqref{E3.11},
$\omega(E_n^{(p)})=-1$. Similar  to the previous cases, based on  the parity conjecture one can claim that  the rank is positive. We know that
these primes are in the form of $8k+1$ and by Dirichlet's theorem,
 there are infinitely many primes of this form. Consequently, the  first part of  Silverman's conjecture is true.\\

\item[(ii)]

 As we mentioned in the proof of  (i), one can easily check that $u$ and $v$ are both odd and then
  $n=u^4+v^4\equiv2\pmod{16}$. Now, we consider two  cases:\\

(a): $p\equiv3\pmod{4}$ and $p^2\equiv1\pmod{16}$. \
\newline
In this case, $\left(\frac{-1}{p}\right)=-1$ and \eqref{E3.1} shows
that  $\epsilon(p^2n)=1$  which together with  \eqref{E3.11}, implies
$\omega(E_n^{(p)})=1$.
 So, under the assumption of the parity conjecture  the rank of $E^{(p)}(\Bbb{Q})$ is an even number.\\

(b): $p\equiv3\pmod{4}$ and $p^2\equiv9\pmod{16}$.\\
\
\newline
Finally in the latest case, we have $\left(\frac{-1}{p}\right)=-1$,
$\epsilon(p^2n)=1$. Therefore, \mbox{$\omega(E_n^{p})=1$}, then the
rank of $E_n^{(p)}(\Bbb{Q})$ in this case  must be even ( by the parity conjecture). \hspace*{6cm}$\Box$
\end{itemize}
The following examples show computations for some specific values
for $n$ and $p$. These computations were done by SAGE \cite{sage}
and Cremona's MWRANK \cite{mwrank} softwares.
 \begin{rem}
Determining the rank of an elliptic curve is a challenging problem and in a lot of cases MWRANK program can not compute the rank, in these cases it gives the upper and lower bounds for the rank.
 \end{rem}
\begin{example}
By taking $n=635318657$  and $p<1000000$ in  Theorem \ref{t1}
\text(i), the maximal rank that we found, was $5$ and there are
exactly two such elliptic curves, namely:
\begin{equation*}
\begin{array}{l}
y^2=x^3-886117685355977x,\\
\\
y^2=x^3-1608116501388523697x.
\end{array}
\end{equation*}
The former curve corresponds to $p=1181$, where $p\equiv1\pmod{4}$
and  $p^2\equiv9\pmod{16}$.
The latter one corresponds to $p=50311$, where $p\equiv3\pmod{16}$
and $p^2\equiv1\pmod{16}$.
\end{example}
\begin{example}
 In Tables  1 and 2,  we summarized the results for  the number $n=635318657$  in  Theorem \ref{t1} \text{(ii)}.
In the former case all primes $p\equiv1\pmod{8}$  in the range
$1000$  have been considered for which  there are no curves with rank $0$.  Table 2 shows the results for   primes $p\equiv3\pmod{8}$ in the range $1000$. In this case we found only $7$ curves having rank  $0$.
\\
\begin{center}
\begin{table}[h]
\caption{Primes $\equiv1\pmod{8}$ and $<1000$}\label{table1}
\renewcommand\arraystretch{1.5}
\noindent\[
\begin{array}{|c|c|}
\hline p& rank\\\hline 89,\  929,\  977&0\leq rank\leq2\\\hline
137,\ 881& rank=2\\\hline
97,\  113,\  233,\  257,\  337,\  353,\  409,&2\leq rank\leq4\\
 433,\  521,\ 601,\  641,\  673,\  761,\  937,\  953&\\\hline
 17,\ 41,\ 73,\ 193,\ 241,\ 281,\ 313,\ 401,\ 449,\ 457,\ 569,&0\leq rank\leq4\\
 577,\ 593,\  617,\ 769,\ 809,\ 857,\ &\\\hline
\end{array}
\]
\end{table}

\begin{table}[h]
\caption{Primes $\equiv3\pmod{8}$ and $<1000$}\label{table2}
\renewcommand\arraystretch{1.5}
\noindent\[
\begin{array}{|c|c|}
\hline p&rank\\\hline 19,\  59,\  179,\  491,\  523,\  587,\
971&rank=0\\\hline
11,\ 83,\  211,\ 251,\  283,\ 331,\  379,\  499,\  547,\  619,& rank=2\\
659,\ 683,\ 811&\\\hline
3,\ 43,\  67,\  107,\  131,\  139,\  163,\  227,\ &\\
307,\  347,\ 419,\  443,\  43,\  467,\  563,\  571,\ 643&0\leq rank\leq2\\
691,\  739,\  787,\  827,\  859,\  883,\  907,\  947&\\\hline
 \end{array}
  \]
\end{table}
\end{center}
 \end{example}
\begin{example}
Let $n=156700232476402$.  We considered primes
\mbox{$p\equiv3\pmod{4}$} in the range the $1000$ and among them,
for
\begin{equation*}
p=59,\ 131,\ 163,\ 211,\ 307,\ 331,\ 347,\ 379,\ 571,\ 587,\ 647,\
691,\ 911,
\end{equation*}
 the rank of $E_n^{(p)}$ is  $0$.
\end{example}
\vskip 2cm



\bigskip
\bigskip


{\footnotesize \pn{\bf Farzali Izadi}\; \\ {Department of
Mathematics}, {Azarbaijan Shahid Madani University,
} {Tabriz 53751-71379, Iran}\\
{\tt Email: farzali.izadi@azaruniv.edu}\\

{\footnotesize \pn{\bf Kamran Nabardi}\; \\ {Department of
Mathematics}, {Azarbaijan Shahid Madani University,
} {Tabriz 53751-71379, Iran}\\
{\tt Email: nabardi@azaruniv.edu}\\
\end{document}